\documentclass[a4paper,11pt]{article}
\usepackage{amssymb,amsmath,amsthm}
\usepackage[dvips]{graphicx}

\newtheorem{thm}{Theorem}[section]
\newtheorem{prop}[thm]{Proposition}

\newtheorem{defin}[thm]{Definition}

\theoremstyle{remark}
\newtheorem{rem}[thm]{Remark}
\newtheorem{ex}[thm]{Example}

\begin{document}

\title{Yang--Baxter maps associated to elliptic curves}
\author{Vassilios G. Papageorgiou{$^{1,}$\footnote{\tt vassilis@math.upatras.gr}} $\phantom{.}$
and Anastasios G. Tongas{$^{2,}$\footnote{\tt atongas@tem.uoc.gr}} \\ \\{\em
$^{1}$ Department of Mathematics, University of Patras, 265 00 Patras, Greece } \\
{\em $^{2}$ Department of Applied Mathematics, University of Crete, 714 09 Heraklion, Greece}}
\maketitle
\begin{abstract}
We present Yang--Baxter maps associated to elliptic curves. They are related to discrete 
versions of the Krichever-Novikov and the Landau-Lifshits equations. 
A lifting of scalar integrable quad--graph equations to two--field equations is also shown.
\end{abstract}

\section{Introduction}

The importance of the Yang-Baxter (YB) equation to a variety of branches 
in physics and mathematics is well known. Its solutions are intimately related to
exactly solvable statistical mechanical models, link polynomials in knot theory,
quantum and classical integrable models, conformal field theories, 
representations of groups and algebras, quantum groups and many others.
More interestingly, YB equation provides various connections among the 
aforementioned disciplines.

Historically, YB equation has its roots in the theory of exactly solvable models in 
statistical mechanics \cite{yang, baxter} and the quantum inverse scattering method \cite{TF}.
For an extensive account of early work on the YB equation see \cite{jimbo}.
In its original form the quantum YB equation is the relation
\begin{equation}
R^{(2,3)}\, R^{(1,3)}\, R^{(1,2)}\, = R^{(1,2)}\, R^{(1,3)}\, R^{(2,3)}\,, 
\label{eq:YBrelV} 
\end{equation}
in $End(V^{\otimes 3})$, for a $k$-linear operator $R:V^{\otimes 2} \mapsto V^{\otimes 2}$, 
where $V$ is a vector space over a field $k$. 
Here, $R^{(1,3)}$ is meant as $R$ 
acting on the first and third factors of the tensor product $V^{\otimes 3}$ and as 
identity on the second, and similarly for $R^{(1,2)}$ and $R^{(2,3)}$. 
Drinfel'd suggested to study the simplest possible solutions of the 
YB equation by replacing $(k\mbox{--Vect},\otimes)$ with $(\mbox{Set},\times)$,
where the YB equation is regarded as an equality of maps in $\mathbb{X}^{3}$
for a finite set $\mathbb{X}$. As it was pointed out in \cite{drin}, this setting 
provides potentially new interesting solutions of the original YB equation
by considering the free module generated by the set $\mathbb{X}$.
Various interesting examples of YB maps such as those arising from geometric crystalls
\cite{eting2}, have revealed a richer structure of the underlying set; $\mathbb{X}$ is an algebraic variety and $R$
is a birational isomorphism. As in \cite{ves1} we refer to 
solutions of the YB equation in $(\mbox{Set},\times)$ simply as YB maps. 
One has to have in mind though that YB maps through a natural 
coupling can be regarded as equations on the edges of graphs.

One of the most distinguished properties of integrable partial differential 
equations is their invariance under Darboux and B\"acklund transformations 
\cite{MaSa}, \cite{RoSch}. 
The nonlinear superposition formulae of the solutions generated by the
B\"acklund-Darboux transformations provide natural discrete versions of the
continuous equations. In turn, the key transformation properties of the 
discrete equations are intimately related to the YB property for maps or
its proper generalization in higher dimensions namely the functional 
tetrahedron equation. A prime example of the latter are the star-triangle 
transformations in electric networks and the Ising model and their 
connection with discrete equations in the three-dimensional lattice 
associated with the Kadomtsev--Petviashvili (KP) hierarchy and its 
modifications \cite{Kashaev}.

Recent interest on solutions of the YB equation for maps has appeared in the 
literature. This is mainly due to the development of the dynamical theory of YB 
maps \cite{ves1}, and classification results of YB maps for $\mathbb{X}=\mathbb{CP}^1$,
in connection to analogous results for two--dimensional integrable discrete 
equations on the square lattice \cite{ABS1, ABS2}. 
The latter correspondence was further investigated
in \cite{VTS, VT} by exploiting the local groups of symmetry 
transformations of the discrete equations.
It was demonstrated that an integrable quadrilateral equation with a 
sufficient $r$-parameter 
symmetry group gives rise to a YB map. However, there exist discrete 
equations which do not admit any local symmetry group. 
Such an equation is a discrete version of the Krichever-Novikov (KN) 
equation \cite{KN1, KN2}, introduced by V. E. Adler in \cite{adler}.

A characteristic feature of the discrete KN equation is that the lattice parameters 
lay on an elliptic curve. 
Remarkably, from the very early times exactly solvable two-dimensional models 
appeared in statistical mechanics it was observed that the solution of many of these models ultimately leads to the introduction of 
elliptic functions, such as the eight-vertex model on the square lattice \cite{baxter2}. 
Thus, it would be interesting to investigate
whether there exist also solutions of the YB equation for maps related to elliptic curves.
Already this problem was addressed in \cite{odes} where a theoretical framework was 
introduced for deriving YB maps from
factorization of matrix polynomials and $\theta$--functions. 

The main aim of the present work is to exhibit YB maps with parameters  
living on elliptic curves and which are associated to integrable partial
differential equations (PDE).
In Section \ref{sec:definitions} we present background material 
on the YB maps. In Section \ref{sec:quad}, we present key transformation
properties of integrable lattice equations, encoded into braid type equations, 
and give a brief account on symmetry aspects of discrete integrable
equations and their usage in deriving YB maps. 
Finally, we present a way for deriving a YB map from integrable lattice 
equations of certain type, without using a local symmetry group.
The latter method is applied to generic integrable lattice equations 
associated to elliptic curves 
such as discrete KN equation and discrete Landau--Lifshits equation and the 
results are presented in Sections \ref{sec:KN} and \ref{sec:LL}, respectively. 
The paper concludes in Section \ref{sec:conclusions} with various comments
and perspectives.

\section{Definitions and notation} \label{sec:definitions}

\subsection*{Braid and YB maps}
In the following we use the notation and terminology introduced in \cite{ves1}
(for a recent review see \cite{ves2}).
Let $\mathbb{X}$ be an algebraic variety, and $R:\mathbb{X} \times
\mathbb{X} \rightarrow \mathbb{X} \times
\mathbb{X}$ a birational isomorphism. Let $R^{(i,j)}:\mathbb{X}^{n} \rightarrow
\mathbb{X}^{n}$ denote the map acting
as $R$ on the components $(i,j)$ of the $n$-fold Cartesian product $\mathbb{X}\times\mathbb{X}\times\cdots\times\mathbb{X}$
and as the identity on all others. More explicitly, for $x,y\in
\mathbb{X}$ let us write
\begin{equation}
R(x,y)=\big(f(x,y),g(x,y)\big). \label{eq:YBgen}
\end{equation}
Then, for $n \geq 2$ and $1 \leq i,j\leq n$, $i\neq j$ the map $R^{(i,j)}$ is given by
\begin{equation}
R^{(i,j)}(x_1,\ldots ,x_n) = \left\{ \begin{array}{l}
(x_1,\ldots,x_{i-1},f(x_i,x_j),x_{i+1}, \ldots,
x_{j-1},g(x_i,x_j),x_{j+1},\ldots x_n)\,,\quad i<j, \\ \\
(x_1,\ldots,x_{j-1},g(x_i,x_j),x_{j+1}, \ldots,
x_{i-1},f(x_i,x_j),x_{i+1},\ldots x_n)\,,  \quad i>j\,. 
\end{array} \right. \label{eq:YBdefmap}
\end{equation}
In particular, for $n=2$ we have $R^{(1,2)}=R$ and
$R^{(2,1)}(x,y)=\big(g(y,x),f(y,x)\big)$.  The latter
map is $R$ conjugated by the permutation map $\sigma$, defined by $\sigma(x,y)=(y,x)$, i.e.
\begin{equation}
R^{(2,1)} = \sigma\, R \, \sigma\,.  \label{eq:R21}
\end{equation}

\begin{figure}[h]
 \centering
 \includegraphics[width=10.0cm,bb=0 399 642 679]{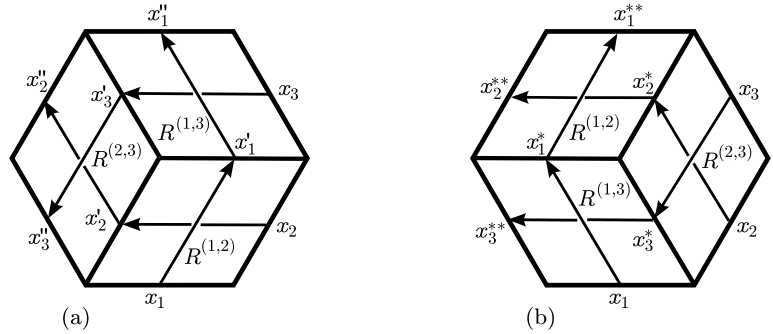} 
\caption{A cubic representation of the YB relation}
\label{fig:3DYB}
\end{figure}

\begin{defin} 
(i) The map $R$ is called a YB map if $R$ satisfies the
YB equation 
\begin{equation}
R^{(2,3)}\, R^{(1,3)}\, R^{(1,2)}\, = R^{(1,2)}\, R^{(1,3)}\, R^{(2,3)}\,, \label{eq:YBrel} 
\end{equation}
regarded as an equality of maps of $\mathbb{X}^{3}$ into itself.

(ii) $R$ is called reversible, or unitary, if it satisfies the condition
\begin{equation}
R^{(2,1)}\,R=\rm{Id}_{\mathbb{X}^{2}}. \label{eq:unit}
\end{equation}

(iii) $R$ is called non-degenerate if the maps from $\mathbb{X}$ into itself defined by
$s \rightarrow f(s,y)$ and $t \rightarrow g(x,t)$  are bijective rational maps 
for any fixed $y$ and $x$, respectively.
\end{defin}
A schematic representation of the YB equation is given by the two decompositions of 
an elementary $3$-cube as depicted in figure \ref{fig:3DYB}. The composition of maps in the 
LHS and RHS of the YB equation (\ref{eq:YBrel}) are given by 
\begin{align}
{\rm (a)}:& \quad (x_1,x_2,x_3) & 
\stackrel{R^{(1,2)}}{\longrightarrow} &&({x'_1},{x'_2},{x_3}) \,&& 
\stackrel{R^{(1,3)}}{\longrightarrow} &&({{x''_1}},{x'_2},{x'_3})\,\, && 
\stackrel{R^{(2,3)}}{\longrightarrow} &&({{x''_1}},{x''_2},{{x''_3}}) \,\, , \nonumber
\\ \label{eq:3DYBinter} \\
{\rm (b)}:& \quad (x_1,x_2,x_3) & 
\stackrel{R^{(2,3)}}{\longrightarrow} &&({x_1},{x^{\ast}_2},{x^{\ast}_3}) && 
\stackrel{R^{(1,3)}}{\longrightarrow} &&({x^{\ast}_1},{x^{\ast}_2},{x^{\ast\ast}_3}) && 
\stackrel{R^{(1,2)}}{\longrightarrow} &&({x^{\ast\ast}_1},{x^{\ast\ast}_2},{x^{\ast\ast}_3}) \,, \nonumber
\end{align} 
respectively. The YB equation guarantees that the images of $(x_1,x_2,x_3)\in \mathbb{X}^{3}$
under the two composition of maps in (\ref{eq:3DYBinter}) are identical, 
thus the two parts of the $3$-cube can be glued together.

\begin{rem}
Defining $B=\sigma R$ the YB equation (\ref{eq:YBrel}) for $R$ translates to the braid
equation 
\begin{equation}
B^{(2,3)}\, B^{(1,2)}\, B^{(2,3)}\, = B^{(1,2)}\, B^{(2,3)}\, B^{(1,2)}\,, 
\label{eq:braidrel} 
\end{equation}
for $B$. 
The unitarity condition (\ref{eq:unit}), in view of (\ref{eq:R21}), takes the form
\begin{equation}
B^2 = {\rm Id}_{\mathbb{X}^{2}} \,,
\end{equation}
and $B$ is said to be an involution. 
\end{rem}

\subsection*{Lax matrices for YB maps}
Instead of a single map one may consider a whole family of YB maps
parametrized by two continuous parameters $\alpha_1,\alpha_2\in\mathcal{C}$, where $\mathcal{C}$ is an 
algebraic set in $\mathbb{C}^2$.
The YB relation then takes the parameter-dependent form
\begin{equation} 
R^{(2,3)}_{(\alpha_2,\alpha_3)} \, R^{(1,3)}_{(\alpha_1,\alpha_3)}\, R^{(1,2)}_{(\alpha_1,\alpha_2)} = R^{(1,2)}_{(\alpha_1,\alpha_2)} \, R^{(1,3)}_{(\alpha_1,\alpha_3)} \, R^{(2,3)}_{(\alpha_2,\alpha_3)} \,, \label{eq:YBcom} 
\end{equation}
and the unitarity (reversibility) condition becomes
\begin{equation}
R^{(2,1)}_{(\alpha_2,\alpha_1)}\,R^{\phantom{(2,1)}}_{(\alpha_1,\alpha_2)} = 
\rm{Id}_{\mathbb{X}^{2}}\,. \label{eq:unitparam}
\end{equation}
In the following we drop the dependence on the parameters and write just 
$R$ for a two-parameter YB map, since we  always consider maps of this type.

Let $L(x;\alpha,\lambda) \in {\rm{Mat}}(r,\mathbb{X})$ 
be a two-parameter family of $r\times r$ matrices depending on $x\in \mathbb{X}$ and
polynomially/rationally on the coordinates of $\alpha,\lambda\in \mathcal{C}$. 
The following notion of Lax matrix for a YB map was introduced in 
\cite{SV}.

\begin{defin}
(i) $L(x;\alpha,\lambda)$ is called a Lax matrix of the YB map 
$R$, if the relation $R({x_1},{x_2})=(x'_1,x'_2)$ implies that
\begin{equation}
L({x_2};\alpha_2,\lambda)\, L({x_1};\alpha_1,\lambda) =  
L({x'_1};\alpha_1,\lambda)\, L({x'_2};\alpha_2,\lambda)\,,
\label{eq:2-fact}
\end{equation}
for all $\lambda\in \mathcal{C}$.
$L(x;\alpha,\lambda)$ is called a strong Lax matrix of 
$R$, if the converse also holds. \\
(ii) $L(x;\alpha,\lambda)$ satisfies the $n$-factorization property if
the identity
\begin{equation}
L({x'_n};\alpha_n,{\lambda})\,\cdots
L({x'_2};\alpha_{2},{\lambda})\,
L({x'_1};\alpha_1,{\lambda}) \equiv 
L({x}_n;\alpha_n,{\lambda})\, \cdots
L({x}_2;\alpha_{2},{\lambda})\,
L({x}_1;\alpha_1,{\lambda})\,, \label{eq:3-fact}
\end{equation}
over $\mathcal{C}$, implies that ${{x}_i}' = {x}_i$, $i=1,\ldots ,n$. 
\end{defin}
\begin{rem} \label{rem:nfact}
The $2$-factorization property of $L$ corresponds to 
the unitarity property of $R$, while the $3$-factorization property to the YB property.
Indeed, the composition of maps
\begin{equation}
\quad (x_1,x_2) \quad
\stackrel{R^{(1,2)}}{\longrightarrow}\quad ({x'_1},{x'_2}) \quad
\stackrel{R^{(2,1)}}{\longrightarrow}\quad ({x''_1},{x''_2})\,, 
\label{eq:R21R12}
\end{equation}
is represented by the matrix factorization
\begin{equation}
L({x_2};\alpha_2,\lambda)\, L({x_1};\alpha_1,\lambda) =  
L({x'_1};\alpha_1,\lambda)\, L({x'_2};\alpha_2,\lambda) =
L({x''_2};\alpha_2,\lambda)\, L({x''_1};\alpha_1,\lambda)\,. 
\label{eq:R21R12-fact}
\end{equation} 
Evidently, the unitary property of $R$ is equivalent to the 
$2$-factorization property of $L$.
On the other hand, the cubic representation of the YB relation 
(see figure \ref{fig:3DYB}) suggests to consider the product 
\begin{equation*}
L({x}_3;{a}_3,{\lambda}) 
L({x}_2;{a}_2,{\lambda}) 
L({x}_1;{a}_1,{\lambda}) \,. 
\end{equation*} 
It can be factorized in two different ways, according to the composition of maps in 
(\ref{eq:3DYBinter}), i.e.
\begin{eqnarray}
\mbox{(a)}\qquad L({x}_3;{a}_3,{\lambda})
L({x}_2;{a}_2,{\lambda}) 
L({x}_1;{a}_1,{\lambda}) &=& 
L({x}_3;{a}_3,{\lambda})
L({x'_1};{a}_1,{\lambda}) 
L({x'_2};{a}_2,{\lambda}) \nonumber \qquad\\ &=&
L({x''_1};{a}_1,{\lambda})
L({x'_3};{a}_3,{\lambda})
L({x'_2};{a}_2,{\lambda}) \nonumber \\ &=&
L({x''_1};{a}_1,{\lambda})
L({x''_2};{a}_2,{\lambda})
L({x''_3};{a}_3,{\lambda})\,, \nonumber
\end{eqnarray}
\begin{eqnarray}
\mbox{(b)}\qquad  L({x}_3;{a}_3,{\lambda})
L({x}_2;{a}_2,{\lambda}) 
L({x}_1;{a}_1,{\lambda}) &=& 
L({x^{\ast}_2};{a}_2,{\lambda})
L({x^{\ast}_3};{a}_3,{\lambda})
L({x}_1;{a}_1,{\lambda}) \nonumber \\ &=&
L({x^{\ast}_2};{a}_2,{\lambda})
L({x^{\ast}_1};{a}_1,{\lambda})
L({x^{\ast\ast}_3};{a}_3,{\lambda}) \nonumber \\ &=&
L({x^{\ast\ast}_1};{a}_1,{\lambda})
L({x^{\ast\ast}_2};{a}_2,{\lambda})
L({x^{\ast\ast}_3};{a}_3,{\lambda})\,. \nonumber
\end{eqnarray}
where we have used the fact that $L$ is a Lax matrix for $R$ in each face of the cube and
the associativity of matrix multiplication.
Thus, the YB property of $R$ is equivalent to the $3$-factorization property of $L$.
\end{rem}
\begin{ex}{\rm Consider the map $R:\mathbb{CP}^1\times\mathbb{CP}^1 \mapsto \mathbb{CP}^1\times\mathbb{CP}^1$
defined by
\begin{equation} 
R(x,y)=\left(y \, \frac{\alpha_1+x\,y}{\alpha_2+x\,y}\,,\, x \, \frac{\alpha_2+x\,y}{\alpha_1+x\,y}\right)\,,  \label{eq:F3map}
\end{equation} 
introduced in \cite{VT}. The above map admits the Lax matrix
\begin{equation*}
L(x;\alpha,\lambda) =  \left[ \begin{array}{cc} x^2 & \lambda\,x \\ x & \alpha \end{array} \right]\,. 
\end{equation*}
It is straightforward to check that the discrete zero curvature 
equation (\ref{eq:2-fact}), implies the map (\ref{eq:F3map}). Conversely, 
according to \cite{SV} a hint for considering the matrix $L$ of the above form, 
is based on the YB map itself. Indeed, one notices that the second component of 
the map can be 
written as a linear fractional transformation induced by the linear transformation
with matrix $L(x;\alpha_1,\lambda)|_{\lambda=\alpha_2}$ on $\left[y\quad 1\right]^T$.
For the above Lax matrix the $n$-factorization property can be proved as follows.
The null space of the linear transformation 
in the LHS of (\ref{eq:3-fact}), for $\lambda=\alpha_1$, is spanned by the vector
$\left[-\alpha_1 \quad x'_1\right]^T$. Similarly, $\left[-\alpha_1 \quad x_1\right]^T$
spans the null space of the RHS linear transformation.
Because of the identity (\ref{eq:3-fact}), we conclude that $x_1={x'_1}$, 
and the rightmost matrices cancel out. 
Therefore, by induction, $L$ satisfies the $n$-factorization property.
The YB map (\ref{eq:F3map}) is simply
related to the $F_{III}$ map obtained in the recent classification \cite{ABS2} 
for the case $\mathbb{X}=\mathbb{CP}^1$.}
\end{ex}

\section{YB maps and integrable quadrilateral equations} \label{sec:quad}

In this section we present first braid transformation properties of 
integrable discrete equations defined on elementary squares. 
Next, we  briefly summarize a method for obtaining YB maps from 
integrable discrete equations on quad-graphs, which is based on the 
existence of a local group of symmetry transformations of the equations. 
Finally, we present another method to the same end which does not 
prerequisite the existence of a local group of symmetry transformations.

\subsection*{Main properties of integrable discrete equations}

\begin{figure}[h]
 \centering
\includegraphics[width=4.0cm,bb=0 0 391 257]{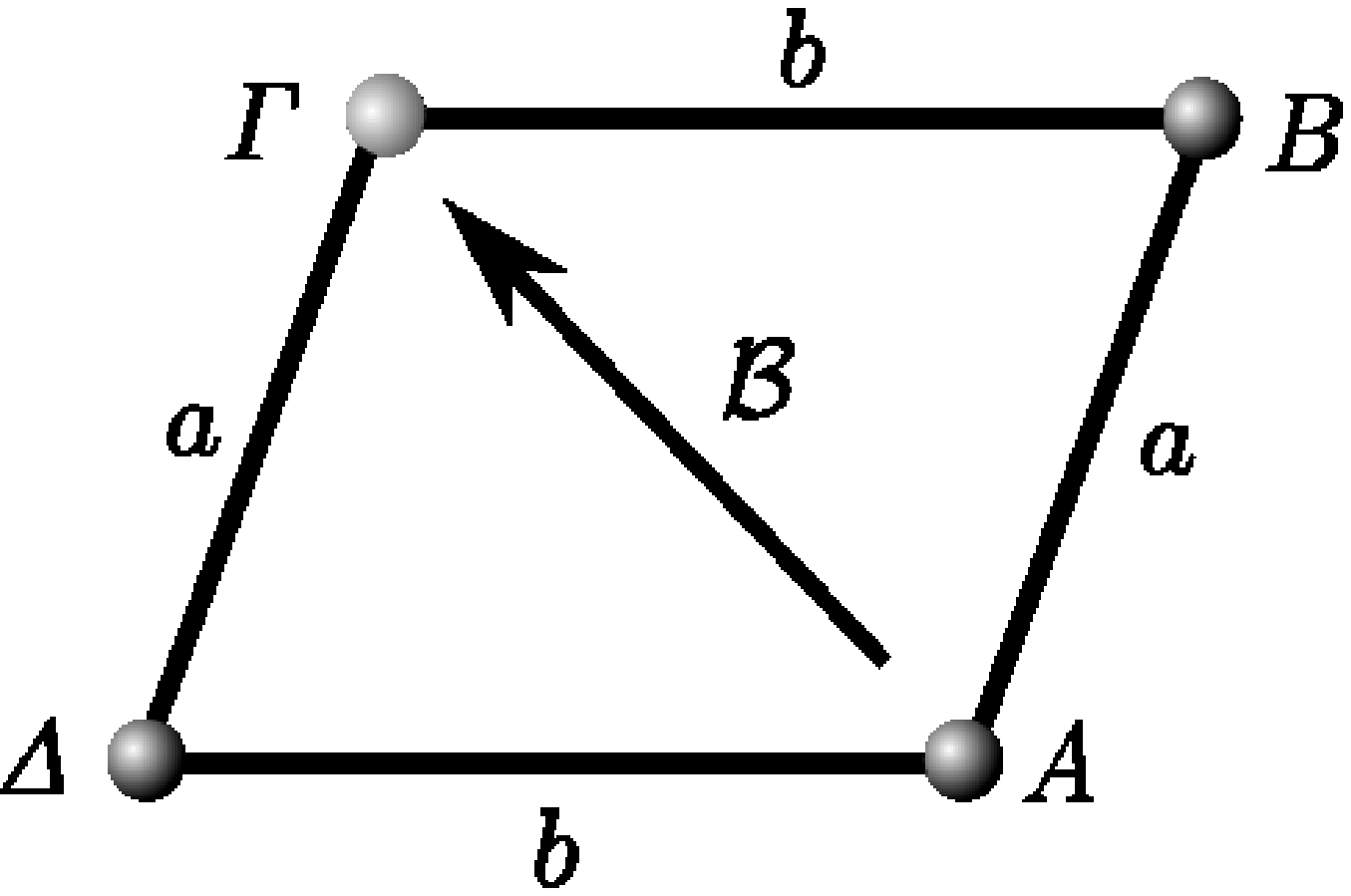}
\caption{An elementary quadrilateral for a quad-graph equation. 
The arrow indicates the flip $(f_{\varDelta},f_{A},f_{B}) \overset{\mathcal{B}}{\rightarrow}
(f_{\varDelta},f_{\varGamma},f_{B})$. }
\label{fig:quad}
\end{figure}

We consider discrete equations on quad-graphs given by 
an algebraic equation 
\begin{equation}
\mathcal{Q}(f_{A},f_{B},f_{\varGamma},f_{\varDelta};\alpha,\beta)=0\,,\label{eq:quad} 
\end{equation}
relating the values of a function $f:\mathbb{Z}^2 \rightarrow \mathbb{X}$ assigned
on the four vertices of an elementary plaquette.
It is assumed that (i) opposite edges on the plaquette carry the same lattice parameter 
$\alpha$, $\beta$ and (ii) equation (\ref{eq:quad})
it can be solved uniquely for each $f_i$, say $f_{\varGamma}$, i.e.
\begin{equation}
f_{\varGamma}=\varphi(f_{A},f_{B},f_{\varDelta};\alpha,\beta)\,.  \label{eq:fgamma}
\end{equation}
In order to make contact with the special properties of the integrable discrete equations 
we interpret equation (\ref{eq:quad}) as a map $\mathcal{B}:\mathbb{X}^{3}\rightarrow\mathbb{X}^{3}$ 
defined by 
\begin{equation}
\mathcal{B}(f_{\varDelta},f_{A},f_{B})=
(f_{\varDelta},f_{\varGamma},f_{B})\,, \label{eq:quadmap} 
\end{equation}
where $f_{\varGamma}$ is given by (\ref{eq:fgamma}) (see fig. (\ref{fig:quad}) ). 

Let us now define 
$\mathcal{B}_j : \mathbb{X}^{n} \rightarrow \mathbb{X}^{n}$ by
\begin{equation}
\mathcal{B}_j  = {\rm Id}_{\mathbb{X}}\times \cdots \times \mathcal{B} 
\times \cdots \times{\rm Id}_{\mathbb{X}}\,, \label{eq:dermap}
\end{equation}
where $\mathcal{B}$ acts on the $j-1$, $j$ and the $j+1$ factors of $\mathbb{X}^n$ 
with parameters $(a_{j-1},a_{j})$.
The key properties of maps associated to integrable discrete equations on 
quad-graphs are the relations 
\begin{equation}
\mathcal{B}_j^2= {\rm Id}_{\mathbb{X}^{n}}\, , \quad 
(\mathcal{B}_j \, \mathcal{B}_{j+1})^3={\rm Id}_{\mathbb{X}^{n}}\,, \quad 
\mathcal{B}_j \, \mathcal{B}_i=\mathcal{B}_i \mathcal{B}_j\,, \quad |i-j|>1 \,,
\end{equation}
see \cite{AY}. The first one means that each transformation $\mathcal{B}_j$ is an involution. 
The second one, in view of the first, yields the following braid-type relation 
\begin{equation}
\mathcal{B}_{j+1}\,\mathcal{B}_j\,\mathcal{B}_{j+1}=
\mathcal{B}_j\,\mathcal{B}_{j+1}\,\mathcal{B}_j.
\label{eq:Braid}
\end{equation}

\begin{figure}[h]
 \centering
 \includegraphics[width=10.0cm,bb= -43 -25 954 344]{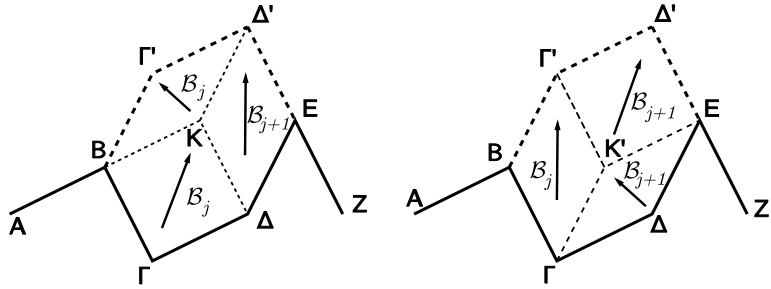} 
\caption{A cubic representation of the braid relation (\ref{eq:Braid}) for quad-graph equations.}
\label{fig:flip}
\end{figure}

The braid relation (\ref{eq:Braid}) guarantees that the three-dimensional 
consistency property \cite{NW}, \cite{BS1},
which nowadays is synonymous with the integrability of a quad-equation holds
(for a recent account on the subject we refer to the monograph \cite{BS2}).

\subsection*{From integrable discrete equations to YB maps via symmetry groups}

Local symmetry groups of transformations of integrable discrete equations 
provide a natural way for obtaining YB maps from them.
The main observation is that the variables of certain YB maps can be chosen 
as invariants of the symmetry group admitted by the corresponding lattice equation.
The symmetry approach was exploited in \cite{VTS}, where it was also shown that 
all classified quatrirational YB maps, for $\mathbb{X}=\mathbb{CP}^1$,
found in \cite{ABS2}, can be constructed from integrable quadrilateral equations. 

\begin{defin} 
Let $G$ be a one-parameter group of transformations on $\mathbb{X}^{3}$, of the form
\begin{equation}
G:\,  (x,y,z) \mapsto \left( X(x;\varepsilon) \,, Y(y;\varepsilon)\,, Z(z;\varepsilon) \right)\,,\qquad 
\varepsilon \in \mathbb{C}\,.
\end{equation}
and $\mathcal{B}$ a map of $\mathbb{X}^{3}$ into itself. 
$G$ is said to be a local (Lie-point) group of symmetry transformations 
of the map $\mathcal{B}$ if $G\circ\mathcal{B}=\mathcal{B} \circ G$, 
for every $\varepsilon \in \mathbb{C}$. 
\end{defin}
Let $\mathbb{X}=\mathbb{CP}^1$ and consider the following map 
\begin{equation}
\mathcal{B}(x,y,z)=
(x,y + \frac{a-b}{x-z},z)\,, \label{eq:quadkdv} 
\end{equation}
which is associated with the discrete KdV equation. 
The corresponding map $\mathcal{B}_j$ defined by
(\ref{eq:dermap}) satisfies the braid type relation (\ref{eq:Braid}). 
Moreover, the map (\ref{eq:quadkdv}) commutes with the group of 
translations given by
\begin{equation}
G_1 : (x,y,z) \mapsto \left(x+ \varepsilon \,, y + \varepsilon\,, z + \varepsilon \right)\,. \label{eq:trans}
\end{equation}
Thus,  $G_1$ is a Lie-point symmetry of the map (\ref{eq:quadkdv}). 
The action of $G_1$ on $\mathbb{X}^{3}$ is regular with one-dimensional orbits, 
thus local coordinates on the set of orbits of $G_1$ are provided by the 
complete set of functionally independent invariants for the group action:
\begin{equation*}
 u=y-x\,,\quad v=z-y\,.
\end{equation*}
Projecting the map (\ref{eq:quadkdv}) to the set of orbits of $G_1$ we obtain the map
\begin{equation}
B(u,v)=\left(u-\frac{a-b}{u+v},v+\frac{a-b}{u+v}\right)\,, \label{eq:kdvbraid} 
\end{equation}
which satisfies the parameter braid relation (\ref{eq:braidrel}). 
Thus, the map
\begin{equation}
R(u,v)=   \sigma \, {B}(u,v) =
\left(v+\frac{a-b}{u+v},u-\frac{a-b}{u+v}\right)\,, \label{eq:kdvYB} 
\end{equation}
is a YB map, known as the Adler map \cite{polygons}. 
The most general local group of symmetry transformations
of the map (\ref{eq:quadkdv}) is $G \cong SO(1,1)$, generated by $G_1$ and the one-parameter subgroups $G_2$, $G_3$ given by the group actions
\begin{eqnarray}
G_2 : (x,y,z) &\mapsto& \left(x - \varepsilon_2 \,, y + \varepsilon_2\,, z - \varepsilon_2 \right)\,, \qquad \nonumber \\
G_3 : (x,y,z) &\mapsto& \left(x \, e^{-\varepsilon_3} \,, y \, e^{\varepsilon_3}\quad\,, 
z \, e^{-\varepsilon_3} \right)\,. \nonumber
\end{eqnarray}
By using similar arguments one may consider the set of orbits of the
subgroups $G_2$, or $G_3$, to obtain other YB maps from the map (\ref{eq:quadkdv}). 
More precisely, we have the following

\begin{prop} \label{prop:symm}
Let $\mathbb{X}=\mathbb{C}^n$ and a map $\mathcal{B}:\mathbb{X}^{3} \rightarrow \mathbb{X}^{3} $
satisfying the braid-type relation (\ref{eq:Braid}). If $\mathcal{B}$  
admits a local group $G$ of symmetry transformations which acts regularly on 
$\mathbb{X}^{3}$ with $n$-dimensional orbits, then the projection
of the map $\mathcal{B}$ to the set of orbits of $G$  
satisfies the braid relation. 
\end{prop}
\begin{proof}
The assumptions for the action of $G$ on $\mathbb{X}^{3}$ guarantee the existence of a
$2n$-dimensional quotient manifold denoted by $\mathbb{X}^{3}/G$,
i.e. the set of all orbits of $G$. Local coordinates $(u,v)\in \mathbb{X}^{2}$
can be chosen by a complete set of functionally independent invariants for the group action, 
see e.g. Theorem 3.18 in \cite{olver}.
 
Let us denote by $B:\mathbb{X}^{2} \mapsto \mathbb{X}^{2} $ 
the projection of the map $\mathcal{B}$ on $\mathbb{X}^{3}/G$. 
The braid property of the map $B$ is inherited by the braid type relation (\ref{eq:Braid}) 
which satisfies the map $\mathcal{B}$. This can be easily deduced from
the cubic representation of the relation (\ref{eq:Braid}) (Figure \ref{fig:flip}). 
It should by noted that the invariants of the group action (YB or braid variables) 
can be naturally assigned to the edges of the 
elementary squares instead of the vertices where the variables of the original map 
$\mathcal{B}$ are assigned to.
\end{proof}
\subsection*{Two--field integrable discrete equations as YB maps} 

The existence of a symmetry group of the map $\mathcal{B}$ provides
us a way to obtain a YB map by using as YB variables the invariants 
of the group action which are naturally attached on the edges of the squares.
As it was shown recently \cite{rahy}, \cite{TTX}, generic integrable 
quad-graph equations, such as the KN discrete equation, do not admit 
any local symmetry group of transformations.
Thus the question arises whether such equations are related to YB maps, as well.
This question is answered in the affirmative in section \ref{sec:KN}.
Proposition \ref{prop:lift} below shows how to cast two--field 
quad--graph equations  of a certain type into YB map form.
Moreover, it motivates a way of lifting an integrable scalar quad--graph equation 
to a two--field one and consequently to recast the equation into a YB map.

Specifically, we consider lattice equations where at each vertex there is a 
two-field $(u,v)\in \mathbb{X}^{2}$ and the defining relations on the quadrilateral 
are (see figure \ref{fig:quad})
\begin{equation}
(u_{\varGamma}\,,\, v_{\varGamma})=\big(F_1(u_A,u_{B},v_{\varDelta};a,b) \,,\,
 F_2(v_A,u_{B},v_{\varDelta};a,b)\big)\,,  \label{eq:AYeqs}
\end{equation}
where $F_1,F_2$ take values in $\mathbb{X}$. This scheme of two-field 
quad-graph equations, although not the generic one 
since it does not involve
all eight values of the fields, arises in the 
superposition formulae of B\"acklund transformations
for two-field integrable PDEs e.g. the nonlinear Schr\"odinger 
system \cite{kono}, \cite{AY}. 
\begin{figure}[h]
 \centering
 \includegraphics[width=5.0cm,bb=108 462 479 680]{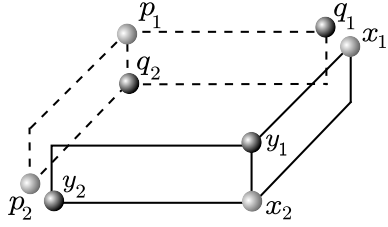} 
 \caption{The quadrilateral with ``thickened" edges as a 
parallelepiped and the assignment of the YB variables for the map (\ref{eq:YBmaplift}).
The variables $x_1$ and $q_1$ are identified, respectively $y_2$ and $p_2$, (6--point scheme).  }
 \label{fig:compina}
\end{figure}
The aim is to recast discrete equations of the form  (\ref{eq:AYeqs}) 
into a YB map form. 
To this end we group the fields appearing in the RHS of 
equation (\ref{eq:AYeqs}) as follows
\begin{eqnarray}
(x_1,x_2) = (u_B,v_A)\,,\qquad  (y_1,y_2) = (u_A,v_\varDelta)\,.
\end{eqnarray}
A pictorial representation of the assignment of the YB variables for the map 
\begin{equation}
 R: \big((x_1,x_2),(y_1,y_2)\big) \rightarrow \big((p_1,p_2),(q_1,q_2)\big) \,,
\end{equation}
(equation (\ref{eq:YBmaplift}) below) is shown in Figure 
\ref{fig:compina} 
where 
\begin{eqnarray}
(p_1,p_2) = (u_\varGamma,v_\varDelta)\,,\qquad  (q_1,q_2) = (u_B,v_\varGamma)\,.
\end{eqnarray}
\begin{prop} \label{prop:lift}
Let the map $\mathcal{B}:\mathbb{X}^{2} \times \mathbb{X}^{2} \times \mathbb{X}^{2} \rightarrow \mathbb{X}^{2} \times \mathbb{X}^{2} \times \mathbb{X}^{2}$, defined by
\begin{equation}
\mathcal{B}\big( (u_{\varDelta},v_{\varDelta}),(u_A,v_A),(u_B,v_B)\big) = 
\big((u_{\varDelta},v_{\varDelta}),(u_\varGamma,v_\varGamma),(u_B,v_B) \big) 
\,, \label{eq:quadmap2} 
\end{equation}
where $(u_\varGamma,v_\varGamma)$ is given by (\ref{eq:AYeqs}), 
satisfies the braid relation (\ref{eq:Braid}). Then the map 
$R: \mathbb{X}^{2} \times \mathbb{X}^{2} \rightarrow \mathbb{X}^{2} \times \mathbb{X}^{2}$ 
defined by
\begin{equation}
R\big((x_1,x_2),(y_1,y_2)\big) = \big( (F_1(y_1,x_1,y_2;\alpha,\beta)\,,\,y_2)\,,(\,x_1\,,\,F_2(x_2,x_1,y_2;\alpha,\beta)\big))\,, \label{eq:YBmaplift}
\end{equation}
satisfies the YB relation.
\end{prop}
\begin{proof}
By straightforward calculations, we derive first the relations for the functions $F_1$, $F_2$
such that $\mathcal{B}$ satisfies the braid type relation. The values $(u_{\varGamma'},v_{\varGamma'})$ 
and $(u_{\varDelta'},v_{\varDelta'})$
are found in two different ways, according to the left and right hand side of the braid relation and are given by
\begin{displaymath}
\begin{array}{lll}
u_K= F_1(u_\varGamma,u_\varDelta,v_B;b,a) &\quad u_{\varDelta'}=F_1(u_\varDelta,u_E,v_K;c,a) &\quad u_{\varGamma'}=F_1(u_K,u_{\varDelta'},v_B;c,b) \\
v_K= F_2(v_\varGamma,u_\varDelta,v_B;b,a) &\quad v_{\varDelta'}=F_2(v_\varDelta,u_E,v_K;c,a) &\quad v_{\varGamma'}=F_2(v_K,u_{\varDelta'},v_B;c,b)
\end{array}
\end{displaymath}
\begin{displaymath}
\begin{array}{lll}
u_{K'}=F_1(u_\varDelta,u_E,v_\varGamma;c,b) &\quad u_{\varGamma'}=F_1(u_{\varGamma},u_{K'},v_B;c,a) &\quad u_{\varDelta'}=F_1(u_{K'},u_E,v_{\varGamma'};b,a) \\ 
v_{K'}=F_2(v_\varDelta,u_E,v_\varGamma;c,b) &\quad v_{\varGamma'}=F_2(v_{\varGamma},u_{K'},v_B;c,a) &\quad v_{\varDelta'}=F_2(v_{K'},u_E,v_{\varGamma'};b,a)
\end{array}
\end{displaymath}
respectively. Thus, we have the following functional relations satisfied by $F_1$, $F_2$:
\begin{eqnarray}
F_1(F_1(u_{\Gamma},u_\varDelta,v_B;b,a),F_1(u_\varDelta,u_E,F_2(v_\varGamma,u_\varDelta,v_B;b,a);c,a),v_B;c,b) = \nonumber \\
F_1(u_\varGamma,F_1(u_\varDelta,u_E,v_\varGamma;c,b),v_B;c,a)\,, 
\label{eq:thmbraid1} \\ \nonumber \\
F_2(F_2(v_\varGamma,u_\varDelta,v_B;b,a),F_1(u_\varDelta,u_E,F_2(v_\varGamma,u_\varDelta,v_B;b,a);c,a),v_B;c,b) = \nonumber \\
F_2(v_\varGamma,F_1(u_\varDelta,u_E,v_\varGamma;c,b),v_B;c,a)\,, 
\label{eq:thmbraid2} \\ \nonumber \\
F_1(u_\varDelta,u_E,F_2(v_\varGamma,u_\varDelta,v_B;b,a);c,a) = \nonumber \\
F_1(F_1(u_\varDelta,u_E,v_\varGamma;c,b),u_E,F_2(v_\varGamma,F_1(u_\varDelta,u_E,v_\varGamma;c,b),v_B,c,a),b,a)\,, \label{eq:thmbraid3}
\\ \nonumber \\
F_2(v_\varDelta,u_E,F_2(v_\varGamma,u_\varDelta,v_B;b,a);c,a)= \nonumber \\
F_2(F_2(v_\varDelta,u_E,v_\varGamma;c,b),u_E,F_2(v_\varGamma,F_1(u_\varDelta,u_E,v_\varGamma;c,b),v_B;c,a);b,a)\,. 
\label{eq:thmbraid4}
\end{eqnarray}
On the other hand, the YB relation for the map (\ref{eq:YBmaplift}) gives the following functional relations for $F_1$, $F_2$
\begin{eqnarray}
F_1(F_1(z_1,y_1,z_2;\beta,\gamma),F_1(y_1,x_1,F_2(y_2,y_1,z_2;\beta,\gamma);\alpha,\gamma),z_2;\alpha,\beta) =
\nonumber \\
F_1(z_1,F_1(y_1,x_1,y_2;\alpha,\beta),z_2;\alpha,\gamma) \,, 
\label{eq:thmYB1} \\ \nonumber \\
F_2(F_2(y_2,y_1,z_2;\beta,\gamma),F_1(y_1,x_1,F_2(y_2,y_1,z_2;\beta,\gamma);\alpha,\gamma),z_2;\alpha,\beta) =
\nonumber \\ 
F_2(y_2,F_1(y_1,x_1,y_2;\alpha,\beta),z_2;\alpha,\gamma)\,, 
\label{eq:thmYB2} \\ \nonumber \\
F_1(y_1,x_1,F_2(y_2,y_1,z_2;\beta,\gamma);\alpha,\gamma) = \nonumber \\
F_1(F_1(y_1,x_1,y_2;\alpha,\beta),x_1,F_2(y_2,F_1(y_1,x_1,y_2;\alpha,\beta),z_2;\alpha,\gamma);\beta,\gamma) \,,
\label{eq:thmYB3} \\ \nonumber \\
F_2(x_2,x_1,F_2(y_2,y_1,z_2;\beta,\gamma);\alpha,\gamma) = \nonumber \\
F_2(F_2(x_2,x_1,y_2;\alpha,\beta),x_1,F_2(y_2,F_1(y_1,x_1,y_2;\alpha,\beta),z_2;\alpha,\gamma);
\beta,\gamma)\,,
\label{eq:thmYB4}
\end{eqnarray}
and two additional equations which are trivially satisfied. 

Making the following substitutions
\begin{displaymath}
\begin{array}{lllll}
&\quad u_\varGamma \mapsto z_1 &\quad u_\varDelta \mapsto y_1 &\quad u_E \mapsto x_1 
&\quad a \mapsto \gamma, \\ 
v_B \mapsto z_2  &\quad v_\varGamma \mapsto y_2 &\quad v_\varDelta \mapsto x_2 &  &\quad b \mapsto \beta \\
 &  & & &\quad c \mapsto \alpha
\end{array}
\end{displaymath}
in equations (\ref{eq:thmbraid1})-(\ref{eq:thmbraid4}), the latter become 
identical to (\ref{eq:thmYB1})-(\ref{eq:thmYB4}), respectively.
\end{proof}
\begin{rem} \label{rem:lift}
Consider the case $F_1=F_2=F$ i.e. 
\begin{equation}
(u_{\varGamma}\,,\, v_{\varGamma})=\big(F(u_A,u_{B},v_{\varDelta};a,b) \,,\,
 F(v_A,u_{B},v_{\varDelta};a,b)\big)\,.  \label{eq:AYeqsF}
\end{equation}
If $u_i=v_i$, $i=A,B,\varDelta$, then from equation (\ref{eq:AYeqsF}) we have
$u_{\varGamma}=v_{\varGamma}$ and the map (\ref{eq:quadmap2}) essentially 
reduces to a single field map, namely 
\begin{equation}
\mathcal{B}^{\downarrow}\big(u_{\varDelta},u_A,u_B\big) = 
\big(u_{\varDelta},u_\varGamma,u_B \big) 
\,, \label{eq:quadmap21} 
\end{equation}
and $\mathcal{B}$ can be thought as a lift of $\mathcal{B}^{\downarrow}$. 
This observation suggests to lift the discrete KN equation to a two-field 
quad-graph equation and then write it as a YB map. 
The lifting process can be applied to all scalar 
integrable quad-equations listed in \cite{ABS1}. 
\end{rem}
\begin{ex}
The simplest equation of the classification in \cite{ABS1} is the discrete 
(potential) KdV equation, namely 
\begin{equation}
f_{n+1,m+1} = f_{n,m} + \frac{a-b}{f_{n+1,m}-f_{n,m+1}} \,, \label{eq:dkdv1}
\end{equation}
where $(n,m)\in\mathbb{Z}^2$. Its lift obtained by equation 
(\ref{eq:AYeqsF}) takes the explicit form
\begin{equation}
u_{n+1,m+1} = u_{n,m} + \frac{a-b}{u_{n+1,m}-v_{n,m+1}} \,, \qquad
v_{n+1,m+1} = v_{n,m} + \frac{a-b}{u_{n+1,m}-v_{n,m+1}} \,, 
\label{eq:dkdv2}
\end{equation}
and satisfies the braid--type relation (\ref{eq:Braid}).
The corresponding YB map obtained by using Proposition \ref{prop:lift} reads
\begin{equation}
(p_1,p_2)=\left(y_1 + \frac{a-b}{x_1-y_2}\,,\,y_2\right)\,,\qquad
(q_1,q_2)=\left(x_1\,,\,x_2 + \frac{a-b}{x_1-y_2}\right) \,. \label{eq:YBkdv}
\end{equation}
The YB map (\ref{eq:YBkdv}) was derived in \cite{koupa} from matrix 
factorization and is symplectic with respect to a canonical structure. 
On the other hand, equations 
(\ref{eq:dkdv2}) are the Euler-Lagrange equations for the discrete 
variational problem associated to the following Lagrangian density 
\begin{equation}
\mathcal{L} =  u_{n+1,m}\,v_{n,m} - u_{n,m}\,v_{n,m+1} + 
(a-b)\, \ln\, (u_{n+1,m} - v_{n,m+1}) \,.
\end{equation}
The problem of the Lagrangian formulation of the quad--graph equations 
classified in \cite{ABS1} has been addressed recently in \cite{loni}.
\end{ex}

\begin{rem}
The YB maps (\ref{eq:YBmaplift}) obtained by this method are non--quadrirational. 
The quadrirationality property of maps as was introduced in \cite{ABS2} 
is equivalent to the nondegeneracy property which is often imposed additionally on the YB maps.
This can be seen immediately since the map $(x_1,x_2) \rightarrow 
(F_1(y_1,x_1,y_2;\alpha,\beta)\,,\,y_2)$ is independent of $x_2$.  
\end{rem}

\section{Lifting discrete KN equation to a YB map} \label{sec:KN}

The master scalar integrable quad equation listed in \cite{ABS1}, in the sense that the rest 
integrable discrete equations can be derived from it by proper degenerations of 
the elliptic curve or limiting procedures, is discrete KN equation \cite{adler}. 
Using the identification $(f_A,f_B,f_{\varGamma},f_{\varDelta})=(x,y,w,z)$ 
on the quadrilateral (Figure \ref{fig:quad}) the latter equation reads
\begin{equation*}
a\,( x y + w z )- b\,( x z + w y)-\frac{a B-b A}{1-a^2 b^2}\big( x w + y z-a b(1+w x y z)\big)=0 \,.
\end{equation*}
This is the form introduced by Hietarinta in \cite{hieta}, where
the parameters $\boldsymbol{a}\equiv(a,A)$ and $\boldsymbol{b}\equiv(b,B)$ lay on 
Jacobi quartics given by
\begin{equation*}
\mathcal{E}=\left\{(a,A)\in \mathbb{C}^2\,:\, {A}^2 = {a}^4 + k \, {a}^2 + 1 \right\} \,,
\end{equation*}
and $k$ is the modulus of $\mathcal{E}$. 
The binary operation $\oplus$ defined by
\begin{equation*}
\boldsymbol{a} \oplus \boldsymbol{b} = \left(\frac{a B+b A}{1-a^2 b^2}\,,\,\frac{(A\,B+k \,a\, b)(1+a^2 b^2)+2 a b (a^2+b^2)}{(1-a^2 b^2)^2}  \right)\,,
\end{equation*}
endows the set $\mathcal{E}$ with an abelian group structure, 
in which $\boldsymbol{e}=(0,1)$ is the identity element and the 
inverse of a point $\boldsymbol{a}=(a,A)$ is the point $-{\boldsymbol a}=(-a,A)$.
In the following we use the notation $\boldsymbol{x}=(x,X)$ etc,
for points in $\mathbb{C}^2$.

\begin{prop} \label{prop:KNlift}
The map 
$R(\boldsymbol{x},\boldsymbol{y})=(\boldsymbol{p},\boldsymbol{q})$ defined by 
\begin{equation}
(p,P)=\left(F(y,x,Y;\boldsymbol{a},\boldsymbol{b}) \,,\,Y\right)\,, \quad
(q,Q)=\left(x \,,\, F(X,x,Y;\boldsymbol{a},\boldsymbol{b}) \right) \,, \label{eq:YBKN}
\end{equation}
where
\begin{equation*}
 F(x,y,z;\boldsymbol{a},\boldsymbol{b}) = 
\frac{(1-a^2 \, b^2) ( b \, z - a \, y)\, x + (a \, B - b \, A) (y \, z - a \, b)}
{(a \, B - b \, A) (a \, b \, y \, z -1)\, x + (1-a^2 \, b^2) ( a \, z - b \, y)}\,,
\end{equation*}
is a unitary YB map, with Lax matrix given by 
\begin{equation}
L(\boldsymbol{x};\boldsymbol{a},\boldsymbol{\lambda}) = 
\rho(\boldsymbol{x};\boldsymbol{a})^{-{1}/{2}} \, W(\boldsymbol{x};\boldsymbol{a},\boldsymbol{\lambda})\,, \label{eq:LaxKN}
\end{equation}
where 
\begin{equation*} W(\boldsymbol{x};\boldsymbol{a},\boldsymbol{\lambda}) = \left[
 \begin{array}{cc} \displaystyle{
-\lambda \, X  - x \,\frac{a \Lambda -A \lambda}{1-a^2 \lambda ^2} }& \displaystyle{
 a\left(1+\lambda \, x \, X\,\frac{a \Lambda -A \lambda}{1-a^2 \lambda ^2} \right) }\\ \\
\displaystyle{-a \left( x \,X+\lambda\,\frac{a \Lambda -A \lambda}{1-a^2 \lambda ^2}\right) }& 
\displaystyle{
\lambda \, x  + X \, \frac{a \Lambda -A \lambda}{1-a^2\lambda ^2} }
\end{array}
\right]\,,
\end{equation*}
and the scalar function $\rho$ is 
\begin{equation*}
 \rho(\boldsymbol{x};\boldsymbol{a})= \left(x^2 \, X^2+1\right) a^2-{x}^2-{X}^2+2\,A\,x\,X\,,
\end{equation*}
$\boldsymbol{a}, \boldsymbol{\lambda} \in \mathcal{E}$.
\end{prop}
\begin{proof}
First we prove that the matrix $L(\boldsymbol{x};\boldsymbol{a},\boldsymbol{\lambda})$
is a Lax matrix for the map (\ref{eq:YBKN}) by showing that the factorization problem
\begin{equation}
L(\boldsymbol{y};\boldsymbol{b},\boldsymbol{\lambda}) \, 
L(\boldsymbol{x};\boldsymbol{a},\boldsymbol{\lambda}) = 
L(\boldsymbol{p};\boldsymbol{a},\boldsymbol{\lambda})\,
L(\boldsymbol{q};\boldsymbol{b},\boldsymbol{\lambda})\,, \label{eq:zerocurv1}
\end{equation} 
is equivalent to equations (\ref{eq:YBKN}). Taking into account that 
$\boldsymbol{\lambda} \in \mathcal{E}$ we have
\begin{equation*}
\big( W(\boldsymbol{y};\boldsymbol{b},\boldsymbol{\lambda})\,
      W(\boldsymbol{x};\boldsymbol{a},\boldsymbol{\lambda}\big)_{ij} = 
\frac{ \sum_{k=0}^{6} \sum_{\ell=0}^{1} S_{ijk\ell}(\boldsymbol{x},\boldsymbol{y};\boldsymbol{a},\boldsymbol{b})\, 
\lambda^k \, \Lambda^\ell}{(a^2 \lambda^2-1) (b^2 \lambda^2-1)}\,.
\end{equation*}
Equating the different powers of $\lambda$, $\Lambda$, 
the matrix equation (\ref{eq:zerocurv1}) is equivalent to the following 
system of algebraic relations
\begin{equation}
\frac{S_{ijk\ell}(\boldsymbol{x},\boldsymbol{y};\boldsymbol{a},\boldsymbol{b})}
{\rho(\boldsymbol{x};\boldsymbol{a})^{1/2} \, \rho(\boldsymbol{y};\boldsymbol{b})^{1/2}} \,
 = 
\frac{S_{ijk\ell}(\boldsymbol{q},\boldsymbol{p};\boldsymbol{b},\boldsymbol{a})}
{\rho(\boldsymbol{p};\boldsymbol{a})^{1/2} \, \rho(\boldsymbol{q};\boldsymbol{b})^{1/2}}\,,
\label{eq:Sijklrat}
\end{equation}
$i,j=1,2$, $k=0,1$ and $\ell=0,\ldots ,6$. We calculate the terms
\begin{align*}
&S_{1100}(\boldsymbol{x},\boldsymbol{y};\boldsymbol{a},\boldsymbol{b}) = a\,b\,x\,(y-X)\,,  \\
&S_{1160}(\boldsymbol{x},\boldsymbol{y};\boldsymbol{a},\boldsymbol{b}) = a^2\,b^2\,Y\,(X-y)\,,  \\
&S_{1201}(\boldsymbol{x},\boldsymbol{y};\boldsymbol{a},\boldsymbol{b}) = a\,b\,(X-y)\,, \\
&S_{1210}(\boldsymbol{x},\boldsymbol{y};\boldsymbol{a},\boldsymbol{b}) = a\,b\,X\,y\, (b\,Y - a\,x)\, + b\,x - a\,Y + a\,B\,y - b\,A\,X\,, \\
&S_{2110}(\boldsymbol{x},\boldsymbol{y};\boldsymbol{a},\boldsymbol{b}) =  
X y (b Y - a x) + a b (b x - a Y) + x Y (a B x - b A y)\,.
\end{align*}
From system (\ref{eq:Sijklrat}) we get
\begin{equation}
\frac{S_{ijk\ell}(\boldsymbol{x},\boldsymbol{y};\boldsymbol{a},\boldsymbol{b})}
{S_{1201}(\boldsymbol{x},\boldsymbol{y};\boldsymbol{a},\boldsymbol{b})} = 
\frac{S_{ijk\ell}(\boldsymbol{q},\boldsymbol{p};\boldsymbol{b},\boldsymbol{a})}
{S_{1201}(\boldsymbol{q},\boldsymbol{p};\boldsymbol{b},\boldsymbol{a}) }\,. \label{eq:Sijkl1201}
\end{equation}
For $(ijk\ell)=(1100)$ and $(ijk\ell)=(1160)$ equations (\ref{eq:Sijkl1201}) 
lead to
\begin{equation}
q=x \,, \qquad P=Y \,, \label{eq:qPeqs}
\end{equation}
respectively and using them, equations (\ref{eq:Sijkl1201}) for 
$(ijk\ell)=(1210),(2110)$ lead to a linear system which is uniquely 
solved for $(p,Q)$, yielding
\begin{equation}
p=F(y,x,Y;\boldsymbol{a},\boldsymbol{b})\,, \qquad 
Q=F(X,x,Y;\boldsymbol{a},\boldsymbol{b}) \,. \label{eq:pQeqs}
\end{equation}
With the solution given by (\ref{eq:qPeqs}), (\ref{eq:pQeqs})
by straightforward calculations we find that system 
(\ref{eq:Sijklrat}) is satisfied. 

Next we prove that the Lax matrix $L$ given by (\ref{eq:LaxKN}) 
satisfies the $n$-factorization property. 
For $\boldsymbol{\lambda} = \boldsymbol{a}$ we obtain
\begin{equation*}
W(\boldsymbol{x};\boldsymbol{a},\boldsymbol{a}) = \alpha 
\left[ \begin{array}{c} 1 \\ {x} \end{array} \right] \left[-{X} \quad 1\right]\,.
\end{equation*}
Thus the kernel of the linear transformation 
\begin{equation*}
L(\boldsymbol{x}_n;\boldsymbol{a}_n,\boldsymbol{a}_1)\cdots L(\boldsymbol{x}_2;\boldsymbol{a}_2,\boldsymbol{a}_1)
\, L(\boldsymbol{x}_1;\boldsymbol{a}_1,\boldsymbol{a}_1)\,,
\end{equation*}
is spanned by the vector
$[1\quad X_1]^T$ which leads us to conclude that ${X_{1}}=X'_{1}$ 
in (\ref{eq:3-fact}).
Likewise, for $\boldsymbol{\lambda}=\boldsymbol{e}$ we have
\begin{equation*}
W(\boldsymbol{x}_n;\boldsymbol{a}_n,\boldsymbol{e})\cdots 
W(\boldsymbol{x}_{2};\boldsymbol{a}_{2},\boldsymbol{e})\,
W(\boldsymbol{x}_1;\boldsymbol{a}_1,\boldsymbol{e}) = \alpha_1 \left( \, 
\prod_{i=2}^n\, \alpha_i \, ({X_{i-1}} - {x_i}) \right) \,  
\left[ \begin{array}{c} 1 \\ {X_n} \end{array} \right] \left[-{x_1} \quad 1 \right]\,.
\end{equation*}
In this case the kernel of the linear transformations in the LHS and RHS of 
(\ref{eq:3-fact}) is spanned by
the vectors $[\,1 \quad x'_1\,]^T$ and $[\,1 \quad {x_1}\,]^T$, respectively.
Thus, $\boldsymbol{x}_1=\boldsymbol{x}'_1$ and consequently
$L(\boldsymbol{x}'_1;\boldsymbol{a}_1,\boldsymbol{\lambda}) = 
L(\boldsymbol{x}_1;\boldsymbol{a}_1,\boldsymbol{\lambda})$ for all 
$\boldsymbol{\lambda}\in \mathcal{E}$. Hence,
the number of matrices in equation (\ref{eq:3-fact}) is reduced by one and
by induction the $n$-factorization property is proved. 
Finally, by remark \ref{rem:nfact} the Proposition is true.
\end{proof}

\section{A discrete Landau-Lifshits equation as a YB map} \label{sec:LL}

In the literature, there exist several discrete versions of the 
Landau-Lifshits equation representing the nonlinear superposition 
formula for the solutions generated by the B\"acklund auto-transformation 
\cite{NiPa}, \cite{AY}, \cite{adler2} . 
Here, we use the one introduced in \cite{AY} and we present the 
end result, namely the corresponding YB map derived from the lattice 
equations by using Proposition \ref{prop:lift}.
The map reads the form 
\begin{equation}
R\big((x,X),(y,Y)\big) = \big( (F_1(y,x,Y;\boldsymbol{a},\boldsymbol{b}),Y)\,,
\,(x,F_2(X,x,Y;\boldsymbol{a},\boldsymbol{b})\big)\,,
\label{eq:YBLL}
\end{equation}
where
\begin{equation*}
F_1(x,y,z;\boldsymbol{a},\boldsymbol{b}) = \frac{K(y,z)\,x - L(y,z)}{M(y,z)\,x+N(y,z)} \,, \quad
F_2(x,y,z;\boldsymbol{a},\boldsymbol{b}) = \frac{K(y,z)\,x + L(y,z)}{-M(y,z)\,x+N(y,z)}\,,
\end{equation*}
\begin{eqnarray*}
K(y,z)-N(y,z) &=& 2\,s_2\, y \, z - (\alpha s_1 + s_3)\,(y-z)-2 \alpha s_2 - 4 \beta s_1\,, \\
K(y,z)+N(y,z) &=& (y+z)\big(a\,b (\alpha s_0 + 3 s_2) + 4 \beta s_1 + 3 \alpha s_2 + s_4\big)/(a-b)\,, \\
L(y,z) &=& s_3 y z + (\alpha s_2 + 2 \beta s_1)(y-z) + 4 \beta s_2 - \alpha^2 s_1\,, \\
M(y,z) &=& s_1 y z + s_2(y-z)-s_3\,,
\end{eqnarray*}
\begin{equation*}
s_m = B \, a^{m-1} + A \, {b}^{m-1} \,. 
\end{equation*}
The parameters $\boldsymbol{a}=(a,A)$, $\boldsymbol{b}=(b,B)$ lay on the 
Weierstra{\ss} elliptic curve
\begin{equation}
{\mathcal{E}}=\left\{(\chi,\mathcal{X})\in \mathbb{C}^2\,:\, \mathcal{X}^2 + {\chi}^3 + \alpha\, \chi + \beta =0 \right\}\,,
\end{equation}
where $\alpha,\beta$ are complex constants, the invariants of the curve.
It should be noted that in contrast to the previous case the functions 
$F_1$, $F_2$ are different reflecting the fact that the discrete equations 
constitute a genuine two-field system. 
The Lax matrix introduced in \cite{AY}, is also a Lax matrix for the map 
(\ref{eq:YBLL}) and is given by
\begin{equation}
L(\boldsymbol{x};\boldsymbol{a},\boldsymbol{\lambda}) = 
\rho(\boldsymbol{x};\boldsymbol{a})^{-1/2} \, W(\boldsymbol{x};\boldsymbol{a},\boldsymbol{\lambda})\,. \label{eq:LaxLL}
\end{equation}
Here,
\begin{equation*}
 \rho(\boldsymbol{x};\boldsymbol{a})= 2\,A (r^2 + x-X + a)\,,
\end{equation*}
the matrix components of $W(\boldsymbol{x};\boldsymbol{a},\boldsymbol{\lambda})$ are given by
\begin{eqnarray*} 
\big(W(\boldsymbol{x};\boldsymbol{a},\boldsymbol{\lambda})\big)_{11} &=&  
(\Lambda+A) r + (\lambda-a)(\lambda + a - X)\,, \\
\big(W(\boldsymbol{x};\boldsymbol{a},\boldsymbol{\lambda})\big)_{12} &=&  
(\Lambda+A)(\lambda + x-X) - (\lambda-a)(\lambda + 2 a) r - 2 A (\lambda - a)\,, \\
\big(W(\boldsymbol{x};\boldsymbol{a},\boldsymbol{\lambda})\big)_{21} &=&  
\Lambda+A - (\lambda-a)r\,, \\
\big(W(\boldsymbol{x};\boldsymbol{a},\boldsymbol{\lambda})\big)_{22} &=&  
-(\Lambda+A) r - (\lambda-a)(\lambda + a + x)\,,
\end{eqnarray*}
and
\begin{equation*}
 r = \frac{1}{2A}\left(x\,X + a(x-X) + \alpha + 2\, a^2\right)\,.
\end{equation*}
The Lax matrix (\ref{eq:LaxLL}) satisfies the $n$-factorization property for $n=2,3$. 
Indeed, first we note that for $\boldsymbol{\lambda}=\boldsymbol{\lambda}_0$, where
\begin{equation*}
\lambda_0 = -2\,a- \mathcal{A}^2 \,, \quad
\Lambda_0 = -A-\mathcal{A}^2(\lambda_0-a) \,, \quad \mathcal{A}=-\frac{3 a^2+\alpha}{2\,A}\,,
\end{equation*}
the matrix $W$ takes the dyadic form
\begin{equation*}
W(\boldsymbol{x};\boldsymbol{a},\boldsymbol{\lambda}_0) = (\lambda_0-a) 
\left[ \begin{array}{c} -\big(2\,A - \mathcal{A} (a-x)\big) \\ \\(a-x)(a+X)\end{array} \right] \,
\Big[ a+X\quad 2\,A- \mathcal{A}(a+X) \Big]\,.
\end{equation*} 
Thus, the kernel of the linear transformation with matrix 
$L(\boldsymbol{x};\boldsymbol{a},\boldsymbol{\lambda}_0)$ is spanned by
the vector $\big[\mathcal{A}(a+X) - 2\,A\quad a+X \big]^T$
from which we conclude that $X=X'$.
Next, for general $\boldsymbol{\lambda}$, the product of two Lax matrices takes the form
\begin{equation*}
\big(L(\boldsymbol{y};\boldsymbol{b},\boldsymbol{\lambda})\,
L(\boldsymbol{x};\boldsymbol{a},\boldsymbol{\lambda}\big)_{ij} = \frac{1}{\rho(\boldsymbol{x})^{1/2} \, \rho(\boldsymbol{y})^{1/2}}
\sum_{k=0}^{3} \sum_{\ell=0}^{1} S_{ijk\ell}(\boldsymbol{x},\boldsymbol{y};\boldsymbol{a},\boldsymbol{b})\, 
\lambda^k \, \Lambda^\ell\,,
\end{equation*}
where we have used that $(\lambda,\Lambda)\in \mathcal{E}$. 
Equating the different powers of $\lambda$, $\Lambda$, 
the matrix equation for the $2$-factorization is equivalent to the 
system of algebraic relations
\begin{equation*}
\frac{S_{ijk\ell}({\boldsymbol{x}'},{\boldsymbol{y}'};
\boldsymbol{a},\boldsymbol{b})}{\rho({\boldsymbol{x}'};\boldsymbol{a})^{1/2} \, \rho({\boldsymbol{y}'};\boldsymbol{b})^{1/2}} = 
\frac{ S_{ijk\ell}(\boldsymbol{x},\boldsymbol{y};\boldsymbol{a},\boldsymbol{b})}
{\rho(\boldsymbol{x};\boldsymbol{a})^{1/2} \, \rho(\boldsymbol{y};\boldsymbol{b})^{1/2}}\,.
\end{equation*}
For $(ijk\ell)=(1221),(1211),(1120)$ we find that
\begin{eqnarray*}
S_{1221}=-(X+y) \,, \qquad  
S_{1211}-S_{1120} = -x (X+y)\,.
\end{eqnarray*}
Hence, $x=x'$ and consequently $L$ satisfies the $2$-factorization property. 

Likewise, the product of three Lax matrices reads
\begin{equation*}
\big(L(\boldsymbol{z};\boldsymbol{c},\boldsymbol{\lambda})\,
L(\boldsymbol{y};\boldsymbol{b},\boldsymbol{\lambda})\,
L(\boldsymbol{x};\boldsymbol{a},\boldsymbol{\lambda}\big)_{ij} = 
\frac{\sum_{k=0}^{4} \sum_{\ell=0}^{1} S_{ijk\ell}(\boldsymbol{x},\boldsymbol{y},\boldsymbol{z};\boldsymbol{a},\boldsymbol{b},\boldsymbol{c})\, 
\lambda^k \, \Lambda^\ell }{\rho(\boldsymbol{x};\boldsymbol{a})^{1/2} 
\, \rho(\boldsymbol{y};\boldsymbol{b})^{1/2} \, \rho(\boldsymbol{z};\boldsymbol{c})^{1/2}}\,.
\end{equation*}
For $(ijk\ell)=(1241),(1221),(1130)$ we find that
\begin{eqnarray*}
S_{1241} = \phantom{x}(X+y)(Y+z)\, , \qquad 
S_{1221}-S_{1130} = x(X+y)(Y+z) \,.
\end{eqnarray*}
From the corresponding ratios we have that $x=x'$, and the number of 
matrices is reduced by one. 
The $3$-factorization property is reduced to the $2$-factorization property, 
which is satisfied, and the map (\ref{eq:YBLL}) is a unitary YB map.

\section{Conclusions} \label{sec:conclusions}

Two families of YB maps with parameters living on elliptic curves are presented. 
Both of them are based on the combinatorics and the geometry of a certain type 
two--field quad--graph system (6-point scheme) that allows to cast them into 
YB map form. 
It is this scheme that suggested the lifting of scalar integrable quad--graph equations to two-field ones and subsequently the derivation of their YB form. 

We end by giving a rough account on YB maps in 
$\mathbb{C}^2\times \mathbb{C}^2$ arising from two--field integrable 
quad--graph equations. 
In \cite{VT} such YB maps were derived 
by exploiting the symmetry groups of the equations
listed in \cite{AY}.
In the present work it is shown how all discrete equations listed 
in \cite{AY} are casted in YB map form (Proposition \ref{prop:lift}). 
This list of YB maps is enhanced by ``lifting"  all integrable 
quad--graph equations listed in \cite{ABS1}, as it was demonstrated
here for the generic equation of the class, namely the discrete KN equation, 
denoted by $Q_4$ in \cite{ABS1}. 
Moreover, the list of YB maps is enriched by considering
the symmetry groups of the lifted discrete equations. 
Thus, it turns out that even in this particular case 
(corresponding to the 6-point scheme) one has already 
quite an amount of YB maps in $\mathbb{C}^2\times \mathbb{C}^2$ 
and the problem of their classification becomes 
interesting in order to: 
i) find the representatives up to equivalence with respect to some 
group of transformations and 
ii) make the list exhaustive.

\subsection*{Acknowledgements}
This work was completed at the Isaac Newton Institute for Mathematical 
Sciences in Cambridge during the programme Discrete Integrable Systems.

\end{document}